\documentclass{amsart}
\newcommand{\note}{\noindent {\bf Notation. }}
\newcommand{\remark}{\noindent {\bf Remark. }}
\newcommand{\corollary}{\noindent {\bf Corollary. }}
\newcommand{\ws}{\hspace{4pt}}
\newtheorem{theorem}{Theorem}
\newtheorem{statement}{Statement}
\newtheorem{proposition}{Proposition}
\newtheorem{lemma}{Lemma}
\newtheorem{defi}{Definition}
\usepackage[]{amssymb,amsopn}
\usepackage[]{amsmath}
\usepackage[]{eucal}
\usepackage[]{latexsym}

\begin{document}

\title[Chromatic Derivatives]{Chromatic Derivatives and Expansions with Weights}
\author{\'A. P. Horv\'ath }

\thanks{Supported by Hungarian National Foundation for Scientific Research, Grant No.  K-100461}

\begin{abstract}   
Chromatic derivatives and series expansions of bandlimited functions have recently been introduced in signal processing and they have been shown to be useful in practical applications. We extend the notion of chromatic derivative using varying weights. When the kernel function of the integral operator is positive, this extension ensures chromatic expansions around every points. Besides old examples, the modified method is demonstrated via some new ones as Walsh-Fourier transform, and Poisson-wavelet transform. Moreover the chromatic expansion of a function in some $L^p$-space is investigated.

\end{abstract} 
\maketitle

\section{Introduction}
A band limited signal, $f(t)$ with finite energy can be represented by different series expansions. As an analytic function, it can be represented by its Taylor series, as an $L^2$-function, it can be expanded to some Fourier series, and by the Whittaker-Shannon-Kotel'nikov sampling theorem it can be represented in the form
$$f(t)\sim \sum_{n=-\infty}^{\infty}f(n)\frac{\sin\pi(t-n)}{\pi(t-n)}.$$
The notion of chromatic derivatives and chromatic expansions have recently been introduced by A. Ignjatovi\'c (cf. e.g. \cite{i0} - \cite{i3}). This alternative representation is possessed of some useful properties as follows. It is of local nature as Taylor series, but it can be handled numerically much better than the Taylor expansion. In some cases it is a Fourier series in some Hilbert space, or its fundamental system of functions forms a Riesz basis or a frame, furtheremore its fundamental system of functions is generated by a single function as Gabor or wavelet systems. The chromatic series expansions of signals associated with Fourier transform have several practical applications, cf. e.g. \cite{ch}, \cite{hb}, \cite{s} and references therein.

The extension of the method to more general integral transforms, differential operators and weighted spaces was introduced by A. Zayed (cf. e.g. \cite{z1} -  \cite{z2}). This type of generalization allows to find connections of the chromatic derivatives and expansions with the Mellin transform technique which was applied to studying the asymptotic behavior of Fourier ransforms of orthogonal polynomials and with some kinds of differential equations cf. \cite{z2} and references therein.

Our aim is to extend chromatic derivatives and chromatic expansions using varying weights. In some cases (for instance Fourier, Walsh-Fourier, Hankel transforms) there is only a single point $x_0$ for wich the kernel function $\Psi(x_0,\cdot)$ is positive, but for instance in cases of Laplace, Bargmann and Poisson-wavelet transforms it has at most finite many sign changes for a large set of $x_0$-s. In these cases we define chromatic series expansions around every $x_0$-s in question. 

The other direction of our investigation is introducing the chromatic expansion of a function $f$ in certain $L^p$-spaces and proving the $L^p$-convergence of its de la Vall\'ee Poussin means. It also turns out, that the rate of convergence of the de la Vall\'ee Poussin means can be estimated by the error of the near-best approximating polynomials of $\hat{f}$ in the weighted space in question.

The article is organized as follows. In Section 2 we give a summary of the method and some results from the literature. By the method of section 4, we prove some convergence theorems in sections 5 and 6. Sections 3 and 6 contain the main curiosity of the paper, namely the modified method, and the examples with varying weights. 

\section{Preliminaries}
In this section we summarize the method and some properties of the chromatic derivative and chromatic expansion. More details can be found in \cite{ch},  \cite{hb}, \cite{i0} - \cite{iz}, \cite{s}, \cite{sw2}, \cite{z1} - \cite{z2} and references therein. Before Theorem 1 we do not take care on function spaces and convergence.

Let $\Psi(x,y)$ be a kernel function and $\mathcal{I}_{\Psi}$ an integral transform on some function spaces, with a measure $\nu$ on a finite or infinite interval $(c,d)$, say 
\begin{equation}f(x):=\mathcal{I}_{\Psi}(\hat{f})(x)=\int_c^d\hat{f}(y)\overline{\Psi(x,y)}d\nu(y).\end{equation}
Let $f(x)$ be defined on $(a,b)$. Let $L$ be a linear differential operator:
\begin{equation}L_x(f)(x)=\sum_{k=0}^md_k(x)\frac{\partial^kf}{\partial x^k}(x),\end{equation}
where $d_k(x)$ are continuous functions on $(a,b)$. Let us suppose that $\Psi(x,y)$ is smooth as a function of $x$, and
\begin{equation}L_x(\overline{\Psi(x,y)})= y\overline{\Psi(x,y)}.\end{equation}
Finally let us assume that there is an $x_0 \in (a,b)$ such that 
\begin{equation}\Psi(x_0,y)\equiv const\neq 0 \ws\ws y\in(c,d),\end{equation} 
(or if $\Psi(x_0,y)\equiv 0$, then let us assume that for some $\alpha$, $\lim_{x\to 0(+)}\frac{\Psi(x,y)}{x^{\alpha}}\equiv const\neq 0$ on $(c,d)$).

(Sometimes it is important that $\mathcal{I}_{\Psi}$ is the inverse transform of an integral transform, that is $\mathcal{I}_{\Psi}=\left(\mathcal{I}_{\Phi}\right)^{-1}$, where $\Phi(x,y)$ is a kernel function on a finite or infinite interval $(a,b)$, and $\mathcal{I}_{\Phi}(f)(y)=:\hat{f}(y)=\int_a^bf(x)\overline{\Phi(x,y)}d\mu(x)$. cf. e.g. \cite{i2} and subsection 3.1 in \cite{z2}.)

\medskip

\noindent {\bf Examples.}
There are several examples in the literature on operators like above (cf. e.g. \cite{i1}, \cite{i2}, \cite{i3}, \cite{z1}, \cite{z2}, \cite{z3}), e. g. the Fourier transform on $\mathbb{R}$: $\overline{\Psi(x,y)}= e^{ixy}$, $L_x=-i\frac{\partial}{\partial x}$, $\Psi(0,y)\equiv 1$; the Laplace transform on $\mathbb{R}_{+}$: $\Psi(x,y)=e^{-xy}$, $L_x=-\frac{\partial}{\partial x}$, $\Psi(0,y)\equiv 1$; the Hankel transform on $\mathbb{R}_{+}$: $\Psi(x,y)=\sqrt{x}J_{\alpha}\left(x\sqrt{y}\right)$, $L_x=-\frac{\partial^2}{\partial x^2}+\frac{\alpha^2-\frac{1}{4}}{x^2}$, $\lim_{x\to 0+}\frac{\Psi(x,y)}{x^{\alpha+\frac{1}{2}}}=\frac{y^{\frac{\alpha}{2}}}{2^{\nu}\Gamma(\alpha+1)}$, that is $d\nu(y)=\frac{dy}{y^{\frac{\alpha}{2}}}$, etc.

\medskip

Let $p(x)$ be an arbitrary polynomial. According to (3) 
\begin{equation}p(L_x)(\overline{\Psi(x,y)})= p(y)\overline{\Psi(x,y)},\end{equation}
and so by (1)
\begin{equation}p(L_x)(f)(x)= \int_c^d\hat{f}(y)p(y)\overline{ \Psi(x,y)}d\nu(y).\end{equation}
Let $w$ be now a positive weight function on $(c,d)$ with finite moments, and let $\{p_n\}_{n=0}^{\infty}$ be the orthogonal polynomials on $(c,d)$ with respect to $w$, that is let
\begin{equation}\int_c^dp_n(y)p_m(y)w(y)d\nu(y)=\delta_{m,n}.\end{equation}
(Let us remark here that usually $d\nu(y)$ is $dy$, except some cases, see for instance tha Hankel transform.)
The $n^{th}$ chromatic derivative of $f$ with respect to $\mathcal{I}_{\Psi}$ and $w$ is
\begin{equation}K^n(f)(x):= \int_c^d\hat{f}(y)p_n(y)\overline{ \Psi(x,y)}d\nu(y),\end{equation}
when the integral exists. Supposing that there is an $x_0 \in (a,b)$ such that $\overline{ \Psi(x_0,y)}\equiv 1$ on $(c,d)$, we have 
\begin{equation}K^n(f)(x_0):= \int_c^d\frac{\hat{f}(y)}{w(y)}p_n(y)w(y)d\nu(y).\end{equation}
Let
\begin{equation}\varphi(x)=\int_c^dw(y)\overline{\Psi(x,y)}d\nu(y),\end{equation}
and
\begin{equation}\varphi_n(x):=K^n(\varphi)(x)=\int_c^dw(y)p_n(y)\overline{\Psi(x,y)}d\nu(y), \ws \ws x \in (a,b).\end{equation}
The chromatic expansion of $f$ is
\begin{equation}f(x)\sim \sum_{n=0}^{\infty}K^n(f)(x_0)K^n(\varphi)(x).\end{equation}
The conditions below ensure the existence of the integrals above, and the convergence of the series. By the previous notations let $(c,d)= (0,d)$, where $0<d\leq \infty$.

\medskip

\begin{theorem}(Th. 5. in \cite{z1})
Let us assume that $\Psi$ and $\hat{f}$ satisfy the following conditions:

a) $\sqrt{w(y)}\Psi(x,y) \in L^2(\mathbb{R}_{+})$

b) $\frac{\hat{f}}{\sqrt{w}} \in L^2(0,d)$

c) If $d=\infty$, we further assume that the integrals
$$\int_0^{\infty}y^n\hat{f}(y)\overline{\Psi(x,y)}dy <\infty, \ws\ws \ws n=0,1,2,\dots$$
   converge uniformly.

Consider $f(x):=\mathcal{I}_{\Psi}(\hat{f})= \int_0^d\hat{f}(y)\overline{\Psi(x,y)}dy$, $x \in (a,b)$, $0<d\leq \infty$.
Then, the chromatic series expansion of $f$ converges pointwise for $x \in (a,b)$, and locally uniformly in $(a,b)$. Moreover
\begin{equation}\Psi(x,y)=\sum_{n=0}^{\infty}\varphi_n(x)p_n(y).\end{equation}\end{theorem}

\medskip

\section{Modification of the method} 

The main point of our modification is the following: instead of (4) we give a weaker assumption. Propery (4) was ensured by the initial condition (3.5) in \cite{z2}, when the chromatic derivatives were asociated with a Sturm-Liouville differential operator. When the chromatic derivatives were asociated with a more general differential operator, the kernel function of the integral had to be one of the fundamental solutions of the initial value problem (cf. \cite{z2} and references therein). By (4), the coefficients of the chromatic expansion of a function are the chromatic derivatives of the function at a certain $x_0$, so the chromatic expansion is located at a certain point. Originally, in the bandlimited case, by the special translation properties of the Fourier transform, the chromatic expansion had the following form
$$f\sim \sum_{n=0}^{\infty}(-1)^nK^n[f](u)K^n[m](t-u),$$
where $u$ is an arbitrary point of the interval, and the convergence of the expansion was independent of $u$ (cf. \cite{i1}). One hand the motivation of the modification is to get expansions around more points. On the other hand, substituting assumption (4) by a weaker one, namely by the following: let $x_0$ be any point for which $\Psi(x_0,y)$ has finite many sign changes, we can extend the notion of  chromatic derivatives and expansions to a wider family of kernel functions (cf. subsection 6.3).

Finally we also modify the notations in connection with the weight, to make it more compatible to the earlier results.

Let $\Psi(x,y)$ be a kernel function which is (piecewise) smooth as a function of $x$ and $\mathcal{I}_{\Psi}$ be an integral transform on a finite or infinite interval  $(c,d)$ such that
\begin{equation}f(x):=\mathcal{I}_{\Psi}(\hat{f})=\int_c^d\hat{f}(y)\overline{\Psi(x,y)}d\nu(y).\end{equation}

Let us assume that the integral above exists and $f(x)$ is defined on $x\in (a,b)$.
Let $L_x$ be a linear differential operator like (2) for which (3) is satisfied. Furthermore let $w$ be a positive weight function, for simplicity with finite moments on $(c,d)$. Let us denote by
\begin{equation} f_w(x):= \mathcal{I}_{\Psi}\left(\hat{f}\sqrt{w}\right).\end{equation}
Let us recall that if $p_n$ is any polynomial of degree $n$, then 

$$p_n(L_x)(f_w)(x)= \int_c^d\hat{f}(y)\sqrt{w(y)}p_n(y)\overline{\Psi(x,y)}d\nu(y).$$

Let $x_0 \in (a,b)$ such that $\overline{ \Psi(x_0,y)}>0$ on  $(c,d)$, let us denote by

\begin{equation}w_{x_0}(y)=w(y) \Psi(x_0,y)^2 .\end{equation}
If $w_{x_0}$ has finite moments, let $p_{x_0,n}=\gamma_{n,x_0}x^n+\dots$ be the $n^{th}$ orthogonal polynomial with respect to $w_{x_0}$, that is
\begin{equation}\int_c^dp_{x_0,n}(y)p_{x_0,m}(y)w_{x_0}d\nu(y)=\delta_{n,m}.\end{equation}
The $n^{th}$ chromatic derivative of $f$ with respect to $x_0$ and $w$ at $x$ is
\begin{equation}K^n_{w_{x_0}}(f)(x):=p_{x_0,n}(L_x)(f_w)(x)= \int_c^d\hat{f}(y)\sqrt{w(y)}p_{x_0,n}(y)\overline{\Psi(x,y)}d\nu(y).\end{equation}
If we write $\hat{f}=g_{x_0}\sqrt{w_{x_0}}$, we have
\begin{equation}K^n_{w_{x_0}}(f)(x_0)=\langle \hat{f},\sqrt{w_{x_0}}p_{x_0,n}\rangle_{d\nu}=\langle g_{x_0}\sqrt{w_{x_0}},p_{x_0,n}\sqrt{w_{x_0}}\rangle_{d\nu},\end{equation}
that is $K^n_{w_{x_0}}(f)(x_0)$ is the $n^{th}$ Fourier coefficient of $g_{x_0}$ with respect to $\{p_{x_0,k}\}_{k=0}^{\infty}$. Usually $p_{x_0,n}\sqrt{w_{x_0}}\in L^{p'}_{d\nu}$ for all $1\leq p'\leq\infty$, $K^n_{w_{x_0}}(f)(x_0)$ exists if $\hat{f}\in L^{p}_{d\nu}$ for some $1\leq p\leq\infty$, and similarly $K^n_{w_{x_0}}(f)(x)$ also exists.
Thus $g_{x_0}$ has the formal Fourier series
$$g_{x_0}\sim \sum_{n=0}^{\infty}K^n_{w_{x_0}}(f)(x_0)p_{x_0,n}.$$
Let
\begin{equation}\varphi_{x_0,n}:=\mathcal{I}_{\Psi}(p_{x_0,n}\sqrt{w_{x_0}}),\end{equation}
which exists usually in the ordinary sense as above, and
\begin{equation}\varphi_{x_0}:=\mathcal{I}_{\Psi}(\Psi(x_0,\cdot))\end{equation}
in ordinary or in distribution sense. Also in ordinary or in distribution sense (cf. \cite{z3}, section 5),
\begin{equation}\varphi_{x_0,n}=K^n_{w_{x_0}}(\varphi_{x_0}).\end{equation}
Thus $f$ has the formal chromatic expansion
\begin{equation}f(x) \sim \sum_{n=0}^{\infty}K^n_{w_{x_0}}(f)(x_0)K^n_{w_{x_0}}(\varphi_{x_0})(x).\end{equation}
With these notations instead of (13) we have, if $\Psi(x,\cdot) \in L^p((c,d))$ with some $1\leq p \leq \infty$, and formally
\begin{equation}\Psi(x,y)\sim \sum_n\varphi_{x_0,n}(x) p_{x_0,n}(y)\sqrt{w_{x_0}(y)}.\end{equation}

When $\Psi(x_0,y)$ has finite many sign changes, according to the sign of $\Psi(x_0,y)$ we decompose $(c,d)$ to finite many subintervals, and proceed as above on every subintervals independently (cf. subsection 6.3 below). 

If $\mathcal{I}_{\Psi}$ is the inverse transform of some integral transform on $(a,b)$ $\mathcal{I}_{\Phi}$, that is $\hat{f}(x)=\int_a^bf(x)\overline{\Phi(x,y)}d\mu(x)$ such that, say
$$\overline{\Phi(x,y)}=\Psi(x,y),$$
and the operator is isometric:
$$\langle \hat{f},\hat{g}\rangle_{d\nu}=\langle f,g\rangle_{d\mu},$$
then by the transform $\mathcal{I}_{\Psi}$ one can get the formal expansion of $f$ with respect to the orthonormal system $\varphi_{x_0,n} :=\mathcal{I}_{\Psi}(p_{x_0,n}\sqrt{w_{x_0}})$. This is the situation for instance in cases of Fourier-, Hankel-, Walsh-, Sturm-Liouville- and  wavelet transforms. 
By (19) it is clear that if the system $\varphi_{x_0,n}$ is complete, the Fourier expansion of an $f\in L^2(a,b)$ with respect to $\varphi_{x_0,n}$ coincides with its chromatic expansion,
$$f(x) = \sum_{n=0}^{\infty}\langle f,\varphi_{x_0,n}\rangle_{d\mu}\varphi_{x_0,n}(x)=\sum_{n=0}^{\infty}K^n_{w_{x_0}}(f)(x_0)K^n_{w_{x_0}}(\varphi)(x).$$
Without an isometry we also have series expansions as above, but $\{\varphi_{x_0,n}\}$ is no longer orthogonal.
In the first case one can use the usual Fourier seies technics to study the convergence properties of the series, in the second case one can take into consideration the properties of the transformation.

It is easy to see that the method cited from the literature and summarized in the previous section is a subcase of this modified one. Examples with varying weights can be found in Section 6.

\medskip

\section{Convergence}

As we have seen (cf. e.g. Th.1 and the cited references) there are $L^2$, pointwise and locally uniform convergence theorems in connection with chromatic series expansions. In this section we give the sketch of thoughts to get $L^p$ and (locally) uniform convergence theorems by de la Vall\'ee Poussin means of the chromatic series, and estimates on the rate of convergence. The theorems and proofs in concrete cases can be found in the next sections.

Let us assume that $\hat{h}$ is a function in some (Banach) space and $\mathcal{I}_{\Psi}(\hat{h})$ exists. Furthermore let us assume that
\begin{equation}\|\mathcal{I}_{\Psi}(\hat{h})\|_{q,\mu}\leq C \|\hat{h}\varrho_1\|_{s,\nu},\end{equation}
for some $1\leq q,s \leq \infty$.
Let $\hat{f}$ be as before, which has  a formal Fourier series with respect to $\{p_{x_0,n}\sqrt{w_{x_0}}\}$. Denoting by $g_{x_0}:=\frac{\hat{f}}{\sqrt{w_{x_0}}}$, $g_{x_0}$ has the expansion $g_{x_0}\sim \sum_{n=0}^{\infty}K^n_{w_{x_0}}(f)(x_0)p_{x_0,n}$. $\hat{S}_n:=\hat{S}_n(g_{x_0},y)$ is the $n^{th}$ partial sum of this series. Let 
$$\hat{\sigma}_n:=\hat{\sigma}_n(g_{x_0},y)=\frac{1}{n}\sum_{k=0}^{n-1}\hat{S}_k, \ws\ws \hat{v}_n:=\hat{v}_n(g_{x_0},y)=2\hat{\sigma}_{2n}-\hat{\sigma}_n$$ 
be the $n^{th}$ Ces\`aro mean and  the $n^{th}$ de la Vall\'ee Poussin mean of this formal series respectively. We chose de la Vall\'ee Poussin mean because of its reproducing property. Let $\varrho_i(y)>0$ ($i=1,2$) and measurable functions on $(c,d)$. Let  us assume further that
\begin{equation}\|\hat{\sigma}_n(g_{x_0})\sqrt{w_{x_0}}\varrho_1\|_{s,\nu} \leq C \|\hat{f}\varrho_2\|_{s,\nu}.\end{equation}
Denoting by 
$$\tilde{v}_n:=\mathcal{I}_{\Psi}(\hat{v}_n(g_{x_0})\sqrt{w_{x_0}}),$$ 
by (25) and (26) we have
$$\|\tilde{v}_n - f\|_{q,\mu} \leq C \|(\hat{v}_n(g_{x_0})-g_{x_0})\sqrt{w_{x_0}}\varrho_1\|_{s,\nu}
$$ $$\leq C\left(\|(\hat{v}_n(g_{x_0}-p_n))\sqrt{w_{x_0}}\varrho_1\|_{s,\nu}+\|(g_{x_0}-p_n)\sqrt{w_{x_0}}\varrho_1\|_{s,\nu}\right)$$ $$\leq C\left(\|(g_{x_0}-p_n)\sqrt{w_{x_0}}\varrho_2\|_{s,\nu}+\|(g_{x_0}-p_n)\sqrt{w_{x_0}}\varrho_1\|_{s,\nu}\right).$$
That is if there is a sequence of polynomials such that the right-hand side of the previous expression tends to zero when $n$ tends to infinity, then the de la Vall\'ee Poussin means of the chromatic expansion tend to the function in the mean.

Similarly
$$|(\tilde{v}_n - f)(x)|\leq \left|\mathcal{I}_{\Psi}((g_{x_0}-p_n)\sqrt{w_{x_0}})(x)\right|+\left|\mathcal{I}_{\Psi}(\hat{v}_n(g_{x_0}-p_n)\sqrt{w_{x_0}})(x)\right|$$
If
\begin{equation}\left|\mathcal{I}_{\Psi}(\hat{h})(x)\right|\leq C(x)\|\hat{h}\varrho_1\|_{s,\nu},\end{equation}
and the polynomials are dense in the space in question, then $v_n$ tends pointwise to $f$. If $C(x)$ has a uniform bound on every compact subintervals of $(a,b)$, then the convergence is locally uniform.

In order to give the estimations in a closed form, let us introduce the following notation:

\medskip

\note Let $w$ be a positive weight function with finite moments on an interval $(a,b)$. Let $fw \in L^p_{(a,b)}$. The error of the best approximating polynomial of degree $n$ of $f$ in the weighted $L^p$-space is
$$E_n(f)_{p,w}=\inf_{p_n\in\mathcal{P}_n}\|(f-p_n)w\|_p.$$
There are several results connected with estimates of $E_n(f)_{p,w}$, cf. e.g. \cite{dl} in the unweighted case, \cite{ms} for Laguerre weights, \cite{dt} for Jacobi weights, \cite{ksz} for Freud weights, \cite{st} for analytic functions, etc.

Collecting the calculations above and supposing that $p_n$ is a near-best approximating polynomial, we can state

\medskip
\begin{theorem} With the notations above, if (25) and (26) fulfil
$$\|\tilde{v}_n - f\|_{q,\mu} \leq C \max_{i=1,2} E_{n}\left(\frac{\hat{f}}{\sqrt{w_{x_0}}}\right)_{s, \sqrt{w_{x_0}}\varrho_i,\nu},$$
and if (27) fulfils
$$|(\tilde{v}_n - f)(x)|\leq C(x),$$
which ensures pointwise, local uniform or uniform convergence according to the properties of $C(x)$.
\end{theorem}

\medskip

\subsection{Lemmas}

In this subsection we collect from the literature all the lemmas which will be useful to prove Theorem 2 in actual weighted spaces. The lemmas are reformulated to our purpose.

Let $w$ be a positive, measurable weight with finite moments on an interval $(c,d)$, $p_{w,n}$ are the orthonormal polynomials on $(c,d)$ with respect to $w$, $\gamma_n>0$ is the leading coefficient of $p_{w,n}$, $\varrho$ is positive and measurable on  $(c,d)$.
$$\Gamma_{n,w}=\max_{k\leq n}\frac{\gamma_{k,w}}{\gamma_{k+1,w}}, \ws\ws\ws \Lambda_{n,w\varrho}=\left\|w\varrho\sum_{k=0}^np_{w,k}^2\right\|_{\infty}.$$
\begin{lemma}(Lemma 1 in \cite{h})
Let us suppose that
\begin{equation}\frac{\Gamma_{n,w}\Lambda_{n,w\varrho}}{n}<K \end{equation}
for some suitable $K$. Let $0<\varrho<M$ a continuous function on $(c,d)$.
If $h\sqrt{\frac{w}{\varrho}} \in L_p$, then for $1<p<\infty$
\begin{equation}\left \| \sigma_n(h)\sqrt{w\varrho}\right\|_p \leq C(p) \|h\sqrt{\frac{w}{\varrho}}\|_p,\end{equation}
\begin{equation}\left \| \sigma_n(h)\sqrt{w\varrho}\right\|_{\infty} \leq C \|h\sqrt{w}\|_{\infty},\end{equation}
\begin{equation}\left \| \sigma_n(h)\sqrt{w}\right\|_1 \leq C \|h\sqrt{\frac{w}{\varrho}}\|_1.\end{equation}

\end{lemma}

\medskip

\begin{lemma}(cf. Theorem 1 in \cite{p}) Let $1\leq p \leq \infty$, $\alpha\geq 0$, $a,b\geq 0$, $W(x)=e^{-\frac{x^2}{2}}x^{\frac{\alpha}{2}}\left(\frac{x}{1+x}\right)^a(1+x)^b$.
$\sigma_n(f,x)$ is the $n^{th}$ $(C,1)$-mean of the Laguerre-$(\alpha)$-Fourier series of a function $f$. Then
$$\|\sigma_n(f,x)W(x)\|_{p,(0,\infty})\leq C \|f(x)W(x)\|_{p,(0,\infty}),$$
where
$$a<\frac{1}{p'}+\min\left(\frac{\alpha}{2},\frac{1}{4}\right);$$
$$b<\frac{1}{p'}+\frac{3}{4}, \ws\ws a+b<2\frac{1}{p'}+\frac{1}{2},\ws\ws \mbox{if} \ws\ws 1\leq p\leq 4,$$
$$b<\frac{1}{3p'}+\frac{5}{4}, \ws\ws a+b<\frac{4}{3p'}+1,\ws\ws \mbox{if} \ws\ws 4< p\leq \infty,$$

\end{lemma}

\corollary 
{\it Let $\alpha > 0$, $w(x)=x^{\alpha}e^{-x}$, $f\sqrt{w} \in L^p$ and let $\sigma_n(f)$ as above. Then
$$\|\sigma_n(f)\sqrt{w}\|_{p}\leq C \|f\sqrt{w}\|_{p}.$$}

\medskip

In the following Lemma $w$ and $v$ are nonnegative functions, the Laplace transform of a function $f$ is defined on $(0,\infty)$ $\mathcal{L}f$ is also defined on $(0,\infty)$, the Hardy operator $Hf$ is given by $Hf(x)=\int_0^xf(t)dt$, and the adjoint operator of $H$ is given by $H^{*}f(x)=\int_x^{\infty}f(t)dt$. Let us denote by $g^{*}$ the decreasing rearrangemant of $g$.

\begin{lemma}(cf. Theorems 1 and 2 of \cite{b})Let $w$ and $v$ be as above. In order for the Laplace transform $\mathcal{L} :L^p(v) \to L^q(w)$ boundedly, set $u(x)=x^{-2}w\left(\frac{1}{x}\right)$.

If  $1<p\leq q<\infty$, it is sufficient that
$$(\mathrm{a}) \ws\ws \left(H\left(\left(\frac{1}{v}\right)^{*\frac{p'}{p}}\right)\right)^{\frac{1}{p'}}\left(H^{*}(u)\right)^{\frac{1}{q}}\leq C,\ws\ws\mbox{or}\ws\ws (\mathrm{b}) \ws\ws \left(H\left(v^{-\frac{p'}{p}}\right)\right)^{\frac{1}{p'}}\left(\mathcal{L}w\right)^{\frac{1}{q}}\leq C,$$
and if  $1<q<p<\infty$ with $\frac{1}{r}=\frac{1}{q}-\frac{1}{p}$, it is sufficient that
$$(\mathrm{c}) \ws\ws I_1=\int_0^{\infty}\left(\frac{1}{v}\right)^{*\frac{p'}{p}}\left(\left(H\left(\left(\frac{1}{v}\right)^{*\frac{p'}{p}}\right)\right)^{\frac{1}{q'}}\left(H^{*}(u)\right)^{\frac{1}{q}}\right)^r\leq \infty,$$
or
$$(\mathrm{d}) \ws\ws I_2=\int_0^{\infty}v^{-\frac{p'}{p}}\left(\left(H\left(v^{-\frac{p'}{p}}\right)\right)^{\frac{1}{q'}}\left(\mathcal{L}w\right)^{\frac{1}{q}}\right)^r\leq \infty.$$
\end{lemma}

\medskip

\section{Examples with a fixed weight}

\subsection{Fourier transform} Chromatic derivatives and expansions with respect to the Fourier transform were investigated by several authors, cf. e.g. \cite{ch}, \cite{hb}, \cite{s}, \cite{i1}, \cite{s}, \cite{z1}. Here we examine the $L^{p'}$-convergence of the de la Vall\'ee Poussin means of the transformed Fourier series. Our starting point is a function $\hat{f} \in L^{p}$ $1\leq p\leq 2$ instead of $L^{2}$, and the rate of convergence is also estimated.

Let us recall that $(c,d)=\mathbb{R}$, $d\nu(y)=dy$  $\Psi(x,y)=\frac{1}{\sqrt{2\pi}}e^{-ixy}$, that is $\Psi(0,y)\equiv \frac{1}{\sqrt{2\pi}}$ thus $w_{x_0}= w_0=\frac{1}{2\pi}w$, $p_{x_0,n}=\sqrt{2\pi}p_{w,n}$. According to (22)
$$\varphi_0(x)=\mathcal{F}^{-1}\left(\frac{1}{\sqrt{2\pi}}\right)=\delta(x).$$
$L_x=- i\frac{\partial}{\partial x}$, $L_x\overline{\Psi(x,y)}=y\overline{\Psi(x,y)}$. 
We have to mention here that it is easy to see that $\varphi_{0,n}=K^n_{w_{0}}(\varphi_{0})$ in distribution sense. Indeed let $\xi$ be infinitely many differentiable and compactly supported on $\mathbb{R}$ and let $p_{0,n}(y)=\sum_{k=0}^na_ky^k$. Then by definition and finally by changing the order of integrals we have 
$$\left(K^n_{w_{0}}(\varphi_{0})\right)\xi=\left(\left(\sqrt{2\pi}\sum_{k=0}^na_k(-i)^k\frac{\partial^k}{x^k}\right)\left(\mathcal{F}^{-1}\left(\frac{\sqrt{w}}{\sqrt{2\pi}}\right)(x)\right)\right)\xi(x)$$ $$=(\sqrt{w})\left(\sum_{k=0}^na_k(i)^k\mathcal{F}^{-1}\left(\xi^{(k)}\right)\right)=(\sqrt{w(y)})\left(\sum_{k=0}^na_ky^k\mathcal{F}^{-1}(\xi)(y)\right)$$ $$=\int_{\mathbb{R}}\sqrt{w(y)}p_{0,n}(y)\frac{1}{\sqrt{2\pi}}\int_{\mathbb{R}}\xi(x)e^{ixy}dxdy=(\varphi_{0,n})\xi.$$
The integral operator has an inverse on $L^2_{\mathbb{R}}$, say, $d\mu(x)=dx$, and Parseval's formula holds true.

\subsubsection{Hermite weight} (Preliminary results can be found e.g. in \cite{sw2}.)

$w(y)=e^{-y^2}$, $p_{0,n}=\sqrt{2\pi}h_n$, where $h_n$-s are the Hermite polynomials.   As it is shown in \cite{h}, in Lemma 1 $\varrho\equiv 1$ with this weight, and all the assumptions are satisfied. That is in (26)  $\varrho_1=\varrho_2\equiv 1$. It is also wellknown that the Fourier operators map $L^p$ to $L^{p'}$, when $1\leq p \leq 2$ and $\frac{1}{p}+\frac{1}{p'}=1$. Moreover the Hermite functions are the eigenfunctions of the integral operator $\mathcal{I}_{\Psi}\left(h_n(y)e^{-\frac{y^2}{2}}\right)=\mathcal{F}^{-1}\left(h_n(y)e^{-\frac{y^2}{2}}\right)=i^nh_n(y)e^{-\frac{y^2}{2}}$. Thus denoting by $g=\frac{\hat{f}}{\sqrt{w}}$, we have $\tilde{v}_n(x)=\mathcal{F}^{-1}(\hat{v}_n(g,y)\sqrt{w(y)})=v_n(f,x)\sqrt{w(x)}$, where $v_n(f,x)$ is the $n^{th}$ de la Vall\'ee Poussin mean of the Hermite expansion of $f$. That is the computations above give back the result
\begin{proposition} Let $1\leq p\leq 2$, $\hat{f} \in L^{p}(\mathbb{R})$. Then recalling the notation $f=\mathcal{F}^{-1}\left(\hat{f}\right)$
\begin{equation} \left\|f-v_n\sqrt{w}\right\|_{p'} \leq C \left\|\hat{f}-\hat{v}_n\sqrt{w}\right\|_{p}\leq C E_{n}\left(\frac{\hat{f}}{\sqrt{w}}\right)_{p, \sqrt{w}},\end{equation}
where $E_{n}(h)_{p, w}$ is the error of the best approximating polynomial of $h$ of degree $n$, in the weigted $L_p$ space. Moreover if $\hat{f} \in L^{1}$, then
\begin{equation}|f(x)-v_n(x)\sqrt{w(x)}|\leq C E_{n}\left(\frac{\hat{f}}{\sqrt{w}}\right)_{1, \sqrt{w}},\end{equation}
which ensures uniform convergence.\end{proposition}

\subsubsection{Freud weights}

First of all let us recall, that $w(x)=e^{-Q(x)}$ is a Freud weight on $\mathbb{R}$, if $Q$ is positive, even and continuous on $\mathbb{R}$, and twice continuously differentiable and $Q^{'}>0$ on $(0,\infty)$ and there are constants $A,B>1$ such that $A\leq \frac{\left(xQ^{'}(x)\right)^{'}}{Q^{'}(x)}\leq B$ on $(0,\infty)$. Since in Lemma 1 $\varrho\equiv 1$ for Freud weights as well (cf. \cite{h}), all the computations above are valid. The only difference is that $p_{w,n}\sqrt{w}$ are not the eigenfunctions of the Fourier transform, so although $\{\varphi_n\}= \{\mathcal{F}^{-1}(p_{w,n}\sqrt{w})\}$ is an orthonormal system, but they are not weighted polynomials with respect to $w$. $\tilde{v}_n$ is the de la Vall\'ee Poussin mean of the expansion of $f$ with respect to $\varphi_n$. Thus by (25), (26) and Lemma 1 
\begin{statement}
If $\hat{f} \in L^{p}(\mathbb{R})$ $1\leq p\leq 2$, then
$$\|f-\tilde{v}_n\|_{p'} \leq C E_{n}\left(\frac{\hat{f}}{\sqrt{w}}\right)_{p, \sqrt{w}},$$ 
and if $\hat{f} \in L^{1}$, then
$$|f(x)-\tilde{v}_n(x)| \leq C E_{n}\left(\frac{\hat{f}}{\sqrt{w}}\right)_{1, \sqrt{w}}.$$
\end{statement}

\subsubsection{Jacobi weights} (Preliminary results can be found e.g. in \cite{ch}, \cite{hb}, \cite{i0}.)

Let us denote by $w(x):=w^{(\alpha,\beta)}(x)=(1-x)^{\alpha}(1+x)^{\beta}$, $\alpha,\beta>-1$.
According to (25) and (26) we have to apply Lemma 1 again. In case of Jacobi weights $\varrho(x)=\sqrt{1-x^2}$, (cf. \cite{h}). 
\begin{statement} Let us assume that $1\leq p \leq 2$, and $\hat{f}, \frac{\hat{f}}{\sqrt{1-x^2}}\in L^p((-1,1))$. Denoting by
$$\tilde{\tilde{v}}_n=\mathcal{F}^{-1}\left(v_n\left(\frac{\hat{f}}{\sqrt{w\varrho}}\right)\sqrt{w\varrho}\right),$$
we have
$$\|\tilde{\tilde{v}}_n-f\|_{p'}\leq CE_n\left(\frac{\hat{f}}{\sqrt{w\varrho}}\right)_{p,\frac{\sqrt{w}}{\sqrt{\varrho}}},$$
and if $\frac{\hat{f}}{\sqrt{\varrho}}\in L^1$ and $\alpha, \beta \geq \frac{1}{2}$,
$$\tilde{v}_n=\mathcal{F}^{-1}\left(v_n\left(\frac{\hat{f}}{\sqrt{w}}\right)\sqrt{w}\right),$$
$$|\tilde{v}_n(x)-f(x)|\leq CE_n\left(\frac{\hat{f}}{\sqrt{w}}\right)_{1,\frac{\sqrt{w}}{\sqrt{\varrho}}}.$$
\end{statement}

\subsection{Walsh-Fourier-Plancherer transform, Laguerre weight}
In this subsection we cite and follow the notations of \cite{sw}. The dyadic field $\mathbb{F}$ is the set of doubly infinite sequences with entries $x_n=0$ or $x_n=1$, and $\lim_{n\to -\infty}x_n=0$, equipped with a suitable sum and product: $x,y\in \mathbb{F}$, $x+y=(|x_n-y_n|)_{ n\in\mathbb{Z}}$, $xy=(z_n)_{ n\in\mathbb{Z}}$, where $z_n=\sum_{i+j=n}x_iy_j \ws (mod\ws 2)$. $(\mathbb{F},+)$ is an abelian group and $(\mathbb{F},\cdot)$ is an abelian semigroup, $(\mathbb{F},+,\cdot)$ is a commutative normed algebra with the norm: $|x|=\sum_{n\in\mathbb{Z}}x_n2^{-n-1}$. Since $(\mathbb{F},+)$ is a locally compact abelian group with compact unite ball (with another norm) (cf. p. 416. in \cite{sw}), there is a unique normalized Haar measure on $\mathbb{F}$. Let $\pi_n(x)=x_n$.  Characters of the additive group are $\Psi_y(x)=(-1)^{\pi_{-1}(xy)}$. The map $x\to|x|$ takes $\mathbb{F}$ onto $[0,\infty)$ and is 1-1 off a countable subset of $\mathbb{F}$. This identification takes the Haar measure to Lebesgue measure on $[0,\infty)$, the characters to generalized Walsh functions (denoting by $\Psi_y(x)$ again), differentiation to dyadic differentiation and induces a dyadic sum and product on $\mathbb{R}_+$ (cf. p. 420. in \cite{sw}). We work on the halfline. If $f\in L^1_{\mathbb{R}_+}$
$$\mathcal{F}f(t)=\hat{f}(t)=\int_0^{\infty}f(y)\Psi_y(t)dy,$$
if $f\in L^2_{\mathbb{R}_+}$, let $f_t(x)=f(x)$ if $x\in [0,t)$, and $f_t(x)=0$ if $x\geq t$. Then $\mathcal{F}f=\lim_{t\to\infty}\hat{f}_t$ in $L^2$-norm, and if  $f=f_1+f_2 \in L^p_{\mathbb{R}_+}\subset L^1_{\mathbb{R}_+}+L^2_{\mathbb{R}_+}$, $1<p< 2$, then $\mathcal{F}f=\hat{f_1}+\mathcal{F}f_2$. With $1\leq p\leq 2$
$$\|\hat{f}\|_{p'}\leq \|f\|_p,$$
cf. Th. 9 on p. 427. in \cite{sw}. Let us recall the notion of dyadic derivative. The dyadic derivative of $f$ at $x$ is 
$$f^{[1]}(x)=\lim_{n\to\infty}{\bf d_n}f(x), \ws\ws \mbox{where}\ws\ws {\bf d_n}f(x)=\sum_{j=0}^{n-1}2^{j-1}\left(f(x)-f(x\dot{+}2^{-j-1}\right),$$
where $\dot{+}$ is the dyadic sum. $f^{[r]}=(f^{[r-1]})^{[1]}$. We have
$$\Psi_y^{[1]}(t)=y\Psi_y(t),$$
cf. p. 421. in \cite{sw}. Thus $(c,d)=\mathbb{R}_{+}$, $d\nu(y)=dy$ $w(y)=y^{\alpha}e^{-y}$ ($\alpha\geq 0$),  $\Psi(x,y)=\Psi_y(x)$, that is $\Psi(0,y)\equiv 1$ thus $w_{0}= w$. 
\begin{equation}K^n_{w_{0}}(\hat{f})(x)= \int_0^{\infty}f(y)y^{\frac{\alpha}{2}}e^{-\frac{y}{2}}l_n^{\alpha}(y)\Psi_y(x)dy,\end{equation}
when $f\sqrt{w}l_n^{\alpha}\in L^1(\mathbb{R}_+)$, and $l_n^{\alpha}(y)$ are the orthonormal Laguerre polynomials. Denoting by $g(y):=f(y)y^{-\frac{\alpha}{2}}e^{\frac{y}{2}}$,
$$K^n_{w_{0}}(\hat{f})(0)=\int_0^{\infty}g(y)y^{\alpha}e^{-y}l_n^{\alpha}(y)dy,$$
wich is the $n^{th}$ Fourier coefficient of $g$ with respect to the Laguerre weight.
$$\varphi_{0,n}(t)=\mathcal{F}(y^{\frac{\alpha}{2}}e^{-\frac{y}{2}}l_n^{\alpha}(y))$$ 
in ordinary sense.
Thus $\hat{f}$ has the chromatic expansion
$$\hat{f}(t)\sim \sum_{n=0}^{\infty}K^n_{w_{0}}(\hat{f})(0)\varphi_{0,n}(t).$$
It is easy to see that
$$\varphi_0(x)=\delta(x), \ws \mbox{and}\ws K^n_{w_{0}}(\varphi_0)=\varphi_{0,n}$$
in distribution sense. Indeed let us define the $L^1$-strong derivative of a function $f$ by ${\bf d}^{[1]}f:= \lim_{n\to\infty}{\bf d_n}f$ in $L^1$ sense. ${\bf d}^{[r]}f={\bf d}^{[1]}\left({\bf d}^{[r-1]}f\right).$ We say that $\xi$ is in the space of test functions, if $\xi$ is compactly supported on $(0,\infty)$ and  $f^{[r]}={\bf d}^{[r]}f$ for all $r\in \mathbb{N}$ and $\xi, \hat{\xi}\in L^1$. (By Egoroff's theorem the second assumption fulfils e.g. if there is a $K$ such that $\left|f^{[r]}(x)\right|<K$ for all $r\in \mathbb{N}$ and $x\in [A,B]:=\mathrm{supp}f\subset(0,\infty)$.) Let $\xi$ be a test function.
$$(\varphi_0) \xi=(\mathcal{F}(1))\xi=\int_0^{\infty}\hat{\xi}(y)\Psi_0(y)dy=\xi(0),$$
where for the last equality Th. 11 on p. 430. in \cite{sw} can be applied since $\xi$ is $W$-continuous at $0$. The second formula can be computed on the same way as before. Thus we have
$$\hat{f}\sim \sum_{n=0}^{\infty}K^n_{w_{0}}(\hat{f})(0) K^n_{w_{0}}(\varphi_0)$$
again. As in e.g. \cite{z2}, (cf. (24) as well) one can derive 
$$\Psi_y(t)\sim \sum_{n=0}^{\infty}\varphi_{0,n}(t)l_n^{\alpha}(y)\sqrt{w(y)}.$$

Let $\tilde{v}_n$ is the $n^{th}$ de la Vall\'ee Poussin mean of the formal series of $\hat{f}$. Then by the Corollary above and by the Hausdorff-Young inequality for Walsh-Fourier transform, we have
\begin{statement}
For all $f\in L^p$, $1\leq p\leq 2$
\begin{equation}\|\hat{f}-\tilde{v}_n\|_{p'}\leq C E_n\left(\frac{f}{\sqrt{w}}\right)_{\sqrt{w},p}.\end{equation}\end{statement}

\medskip

\section{Examples with varying weights}
In this section we examine transforms with positive kernels and one with a kernel with a single sign change. Here the weight function varies with the point $x_0$, and we give the chromatic derivatives and expansions of the functions in every points $x_0$ of the domain with respect to weights depend on the points.

\subsection{Laplace transform, Laguerre-type weigths}

To get a positive weight on the positive axis, the Laplace transform is restricted to $x,y \in \mathbb{R}_+$. Thus $(c,d)=\mathbb{R}_{+}$, $d\nu(y)=dy$ $w(y)=y^{\alpha}e^{-y}$ $(\alpha\geq 0$),  $\Psi(x,y)=e^{-xy}$. That is $\Psi(0,y)\equiv 1$ thus $w_{0}= w$. This case was studied by A. Zayed cf. \cite{z1}. Since $\Psi(x,y)>0$, 
$$w_{x_0}(y)=w(y) \Psi(x_0,y)^2 =y^{\alpha}e^{-y(2x_0+1)},$$
and
$$p_{x_0,n}(y)=(2x_0+1)^{\frac{\alpha+1}{2}}l_n^{\alpha}((2x_0+1)y),$$
where $l_n^{\alpha}(y)$ are the orthonormal Laguerre polynomials again. $L_x=-\frac{\partial}{\partial x}$. Turning back to the traditional notation we write $f$ instead of $\hat{f}$ and $\bar{f}$ instead of $f$
\begin{equation}K^n_{w_{x_0}}(\bar{f})(x_0)= \int_0^{\infty}f(y)y^{\frac{\alpha}{2}}e^{-y(x_0+\frac{1}{2})}(2x_0+1)^{\frac{\alpha+1}{2}}l_n^{\alpha}((2x_0+1)y)dy.\end{equation}
\begin{equation}\varphi_{x_0}(x)=\mathcal{I}_{\Psi}\left(e^{-x_0s}\right)=\frac{1}{x+x_0},\end{equation}
$$\varphi_{x_0,n}(x)=\int_0^{\infty}y^{\frac{1}{2}\alpha}e^{-y(x_0+\frac{1}{2})}(2x_0+1)^{\frac{\alpha+1}{2}}l_n^{\alpha}((2x_0+1)y)e^{-xy}dy=K^n_{w_{x_0}}\left(\varphi_{x_0}\right)(x)$$ $$=p_{x_0,n}(L_x)\left(\varphi_{x_0,w}\right)=p_{x_0,n}(L_x)\left(\mathcal{L}\left(y^{\frac{\alpha}{2}}\right)\left(x+x_0+\frac{1}{2}\right)\right)$$ \begin{equation}=p_{x_0,n}(L_x)\left(\frac{\Gamma\left(\frac{\alpha}{2}+1\right)}{\left(x+x_0+\frac{1}{2}\right)^{\frac{\alpha}{2}+1}}\right),\end{equation}
and
\begin{equation}\bar{f} \sim \sum_{n=0}^{\infty}K^n_{w_{x_0}}(\bar{f})(x_0)K^n_{w_{x_0}}\left(\varphi_{x_0}\right)(x).\end{equation}
$\{\varphi_{x_0,n}(x)\}$'s are no longer orthogonal, but they form a Riesz basis of $L^2_{(0,\infty)}$, cf. e.g. Th. 1 of \cite{b} and Th. 6.2 in \cite{w}. $\tilde{v}_n$ is the de la Vall\'ee Poussin mean of the formal series above. 

By the classical result connected with Laplace transform (cf. e.g. Lemma 3 (a) with $v=w\equiv 1$) and by the previous Corollary we have 
\begin{proposition}
Let $1\leq p \leq 2$. For an $f \in L^p_{(0,\infty)}$, 
$$\|\tilde{v}_n - \bar{f}\|_{p',(0,\infty)} \leq C \|(\hat{v}_n(g_{x_0})-g_{x_0})\sqrt{w_{x_0}}\|_{p,(0,\infty)}\leq C E_n\left(\frac{f(y)}{\sqrt{w_{x_0}}}\right)_{p,\sqrt{w_{x_0}}}.$$\end{proposition}
By Lemma 2 and Lemma 3 (a), (c) the following theorem holds true:
\begin{theorem}
\begin{equation}\left(\int_0^{\infty}|\tilde{v}_n(x) - \bar{f}(x)|^q\varrho(x)dx\right)^{\frac{1}{q}}\leq C E_n\left(\frac{f(y)}{\sqrt{w_{x_0}}}\right)_{p,\sqrt{w_{x_0}}\omega},\end{equation}
where with $0\leq b\leq a<\frac{1}{p'}$,
$$\omega(y)=\left(\frac{(2x_0+1)y}{1+(2x_0+1)y}\right)^a(1+(2x_0+1)y)^b.$$
If $1<p\leq q <\infty$,
$$\varrho(x)=(1+x)^{-(a-b)q}x^{q\left(\frac{1}{p'}-b\right)-1},$$
and if $1<q<p<\infty$, with $0<\delta<\varepsilon$ arbitrary
$$\varrho(x)=(1+x)^{-(a-b)q-\varepsilon}x^{(a-b)q-\frac{aq+1}{p}+\delta}.$$
\end{theorem}

\medskip

\remark
If $1\leq p \leq 2$ and $q=p'$, $a=b=0$, (40) gives back Prop.2.

\medskip

\proof We have to apply (a) and (c) of Lemma 3  with $v=\omega^p$ and $w=\varrho$.
Observing that if $a,b\geq 0$ $\left(\frac{1}{v}\right)^{*}=\frac{1}{v}$, we have
\begin{equation} \left(\frac{1}{v}\right)^{*\frac{p'}{p}}(y)=((2x_0+1)y)^{-ap'}(1+(2x_0+1)y)^{(a-b)p'}.\end{equation}
Thus by $0<b<a<\frac{1}{p'}$
$$ \left(H\left(\left(\frac{1}{v}\right)^{*\frac{p'}{p}}\right)\right)^{\frac{1}{p'}}$$
\begin{equation}=\left(x\int_0^{1}(1+Cxt)^{(a-b)p'}(Cxt)^{-ap'}dt \right)^{\frac{1}{p'}}\leq Cx^{\frac{1}{p'}-a}(1+x)^{a-b}.\end{equation}
Here and further on $C$ may depend on $x_0$, $p$ and on other fixed parametrs, but it is always independent of $n$. If $1< p\leq q < \infty$,
\begin{equation} \left(H^{*}(u)\right)^{\frac{1}{q}}=\left(\int_x^{\infty}(1+t)^{(b-a)q}t^{-1-q\left(\frac{1}{p'}-a\right)}dt\right)^{\frac{1}{q}}\leq C(1+ x)^{-(a-b)}x^{-\left(\frac{1}{p'}-a\right)}.\end{equation}
That is the assumption (a) of Lemma 3 is satisfied. Let us recall that $\hat{\sigma}_n(h)$ be the $n^{th}$ (C,1) mean of the Fourier series of $h$ with respect to $\{p_{x_0,n}\}$ and  $\tilde{\sigma}_n(h)=\mathcal{L}\left(\hat{\sigma}_n(h)\sqrt{w_{x_0}}\right)$. As in Section 4
$$\left(\int_0^{\infty}|\bar{f}-\tilde{v}_n|^q\varrho\right)^{\frac{1}{q}}\leq C  \|(g_{x_0}-\hat{v}_n)\sqrt{w_{x_0}}\omega\|_p$$ $$
\leq C \left(\|(g_{x_0}-p_n)\sqrt{w_{x_0}}\omega\|_p+\|\hat{\sigma}_n(g_{x_0}-p_n)\sqrt{w_{x_0}}\omega\|_p\right).$$
To apply Lemma 2, let us denote by $h_{*}(y)=h\left(\frac{y}{2x_0+1}\right)$ and by $\sigma_n(h)$ be the $n^{th}$ (C,1) mean of the Fourier series of $h_*$ with respect to $\{l_{n}^{\alpha}\}$. By the substitution $(2x_0+1)y=u$ and by Lemma 2 we have
$$\|\hat{\sigma}_n(h)\sqrt{w_{x_0}}\omega\|_p = C\|\sigma_n(h_*)\sqrt{w}\omega_*\|_p \leq C\|h_*\sqrt{w}\omega_*\|_p=C\|h\sqrt{w_{x_0}}\|_{p,\omega^p}.$$
By this observation, finally we have
\begin{equation}\left(\int_0^{\infty}|\bar{f}-\tilde{v}_n|^q\varrho\right)^{\frac{1}{q}}\leq C \|(g_{x_0}-p_n)\sqrt{w_{x_0}}\omega\|_{p}.\end{equation}
If $p_n$ is a near-best approximating polynomial, we get (40).
Next we check assumption (c) of Lemma 3. Similarly to (43)
$$\left(H^{*}(u)\right)^{\frac{p}{p-q}}=\left(\int_x^{\infty}(1+t)^{(b-a)q-\varepsilon}t^{\frac{aq+1}{p}-2+(\varepsilon-\delta)}dt\right)^{\frac{p}{p-q}}$$
\begin{equation} \leq C(1+ x)^{\frac{(-(a-b)-\varepsilon)p}{p-q}}x^{\left(\frac{aq+1}{p}-1+(\varepsilon-\delta)\right)\frac{p}{p-q}}.\end{equation}
By (41),(42) and (45) $I_1$ is convergent. We can finish the proof as previously.

\medskip

By Lemma 2 and Lemma 3 (b), (d) the following theorem holds true:
\begin{theorem}
\begin{equation}\left(\int_0^{\infty}|\tilde{v}_n(x) - \bar{f}(x)|^q\varrho(x)dx\right)^{\frac{1}{q}}\leq C E_n\left(\frac{f(y)}{\sqrt{w_{x_0}}}\right)_{p,\sqrt{w_{x_0}}\omega}.\end{equation}
If $1<p\leq q <\infty$, let $0\leq b= a<\frac{1}{p'}$ and 
$$\omega(y)=\left((2x_0+1)y\right)^a, \ws\ws \mbox{and} \ws\ws \varrho(x)=x^{q\left(\frac{1}{p'}-a\right)-1}.$$
If $1<q<p<\infty$, let $0\leq b, \frac{1}{q'}\leq a<\frac{1}{p'}$, $\omega(y)=\left(\frac{(2x_0+1)y}{1+(2x_0+1)y}\right)^a(1+(2x_0+1)y)^b$, and
$$\varrho(x)=x^{\frac{(a-b)pq-aq-1}{p}-\delta},$$
where $\delta>0$ arbitrary.
\end{theorem}

\medskip

\proof 
If $1<p\leq q <\infty$, by (42) we need that $\mathcal{L}(\varrho)(y)\leq C y^{-q\left(\frac{1}{p'}-a\right)}$, which is satisfied by the assumption. 

If $1<q<p<\infty$, by (41), (42) and the assumptions on $a$, the exponent of $x$ in $I_2$ is greater than $-1$, so we have to deal with the infinite part of the integral. $I_2$ is finite, if  $\mathcal{L}(\varrho)(y)\leq C(c+ y)^{-\beta}$, where $\beta=1+\frac{-1-aq+(a-b)pq}{p}-\delta$, and this is also satisfied by the assumption. The proof can be finished in both cases as earlier.

\medskip

\subsection{Bargmann transform, Hermite-type weights}

The Bargmann transform was introduced in \cite{ba}, and it was examined in point of chromatic derivatives in \cite{z2}. As it was pointed out in \cite{z2} (cf. also the references therein), the Bargmann transform resurfaced again in recent years in connection with Gabor and Zak transforms. As in the previous subsection, to get positive weights, the Bargmann transform of a function is examined as a function on $\mathbb{R}^n$.
\begin{defi} Let $f: \mathbb{R}^n \to \mathbb{C}$. The Bargmann transform $\mathcal{A}(f)$ of $f$ is defined by
$$\mathcal{A}(f)(z)=F(z)=\int_{\mathbb{R}^n}f(\zeta)k(z,\zeta)d\zeta,$$
where $\zeta, z \in \mathbb{R}^n$ and
$$k(z,\zeta):=k_n(z,\zeta)=\frac{1}{\pi^{\frac{n}{4}}}e^{-\frac{\|z\|^2+\|\zeta\|^2}{2}+\sqrt{2}\langle z,\zeta \rangle}.$$\end{defi}
The corresponding differential operator is (cf. \cite{z2})
$$L_i=\frac{1}{\sqrt{2}}\left(\frac{\partial}{\partial z_i}+z_i\right).$$
Let $\alpha$ be a multiindex, $\alpha=k_1,k_2,\dots,k_n$, and $p_{\alpha}(x)=p_{k_1}(x_1)p_{k_2}(x_2)\dots p_{k_n}(x_n)$ is a polynomial with $n$ variables and of degree $|\alpha|=k_1+\dots+k_n$. It is easy to see that
$$L_ik(z,\zeta)=\zeta_ik(z,\zeta).$$
Thus denoting by $p_{\alpha}(L)=p_{k_1}(L_1)p_{k_2}(L_2)\dots p_{k_n}(L_n)$
$$p_{\alpha}(L)k(z,\zeta)=p_{\alpha}(\zeta)k(z,\zeta).$$
Let $w(\zeta)=w_1(\zeta_1)\dots w_n(\zeta_n)$ a weight function, such that $w_i(\zeta_i)$ are positive and have finite moments. Let $z^0 \in \mathbb{R}^n$, and
$$w_{z^0}(\zeta)=w(\zeta)k^2(z^0,\zeta)=\prod_{i=1}^nw_{z^0_i}(\zeta_i)=\prod_{i=1}^n\frac{1}{\sqrt{\pi}}w_i(\zeta_i)e^{-\left((z_i^0)^2+\zeta_i^2\right)+2\sqrt{2}z_i^0\zeta_i}.$$ 
Let $\{p_{z^0_i,m}\}_m$ the orthonormal polynomials with respect to $w_{z^0_i}$, and $p_{z^0,\alpha}(\zeta)=\prod_{i=1}^np_{z^0_i,k_i}(\zeta_i)$. Obviously $\int_{\mathbb{R}^n}p_{z^0,\alpha}(\zeta)p_{z^0,\beta}(\zeta)w_{z^0}(\zeta)d\zeta=\delta_{\alpha,\beta}$.
The $\alpha^{th}$ chromatic derivative of $F$ with respect to $w_{z^0}$ at $z$ is defined as
$$K^{\alpha}_{w_{z^0}}(F)(z):=\int_{\mathbb{R}^n}f(\zeta)\sqrt{w(\zeta)}p_{z^0,\alpha}(\zeta)k(z,\zeta)d\zeta,$$
and
$$K^{\alpha}_{w_{z^0}}(F)(z^0):=\int_{\mathbb{R}^n}f(\zeta)\sqrt{w_{z^0}(\zeta)}p_{z^0,\alpha}(\zeta)d\zeta=c_{\alpha}.$$
Let
$$\varphi_{z^0}(z)=\mathcal{A}(k(z^0,\zeta))= \frac{1}{\sqrt{\pi}^n}e^{-\frac{\|z^0\|^2+\|z\|^2}{2}}\prod_{i=1}^n\int_{\mathbb{R}}e^{-\zeta_i^2+\sqrt{2}(z^0_i+z_i)\zeta_i}d\zeta_i=e^{\langle z^0,z\rangle}.$$
$$\varphi_{z^0,\alpha}(z)=\mathcal{A}(\sqrt{w_{z^0}(\zeta)}p_{z^0,\alpha}(\zeta))(z)=K^{\alpha}_{w_{z^0}}\left(\varphi_{z^0}\right)(z).$$
Then $F$ has the chromatic expansion
\begin{equation}F(z)\sim \sum_{\alpha}c_{\alpha}\varphi_{z^0,\alpha}(z)=\sum_{\alpha}K^{\alpha}_{w_{z^0}}(F)(z^0)K^{\alpha}_{w_{z^0}}\left(\varphi_{z^0}\right)(z).\end{equation}
Since $k(z,\zeta)\in L^2(\mathbb{R}^n)$ we have the Fourier expansion of the kernel as in (24),
\begin{equation}\frac{1}{\pi^{\frac{n}{4}}}e^{-\frac{\|z\|^2+\|\zeta\|^2}{2}+\sqrt{2}\langle z,\zeta \rangle}\sim \sum_{\alpha}\varphi_{z^0,\alpha}(z)p_{z^0,\alpha}(\zeta)\sqrt{w_{z^0}(\zeta)}.\end{equation}

Restricted to $\mathbb{R}^n$, the inner product of the transformed functions is (cf. \cite{z2})
$$\langle F, G\rangle_{\mathcal{F}}=\int_{\mathbb{R}^n}\frac{1}{\pi^\frac{n}{2}}F(z)e^{-\frac{\|z\|^2}{2}}\frac{1}{\pi^\frac{n}{2}}\overline{G(z)}e^{-\frac{\|z\|^2}{2}}d^nz,$$ and the Bargmann transform maps an $f\in L^2\left(\mathbb{R}^n\right)$ to $F(z)$ such that\\ $\|F\|_{\mathcal{F}}=\left\|\frac{1}{\pi^\frac{n}{2}}F(z)e^{-\frac{\|z\|^2}{2}}\right\|_2$ is finite (cf. Th. 8.2 in \cite{z2}). Thus if $f\in L^2(\mathbb{R}^n)$ (47) and (48) converges in $\mathcal{F}$-norm (cf. Th. 9.3 of \cite{z2}). A similar estimation of norms can be given in $L^p$-norm ($z\in \mathbb{R}^n$). Indeed, let us observe that
$$F(z)e^{-\frac{\|z\|^2}{2}}=\int_{\mathbb{R}^n}f(\zeta)\frac{1}{\pi^{\frac{n}{4}}}e^{-\left\|z-\frac{\zeta}{\sqrt{2}}\right\|^2}d\zeta=(f*g)(\sqrt{2}z),$$
where $g(u)=\left(\frac{4}{\pi}\right)^{\frac{n}{4}}e^{-\frac{\|u\|^2}{2}}$. Thus, by Young's inequality and by changing the variables we have
\begin{statement}
$$\|\frac{1}{\pi^\frac{n}{2}}F(z)e^{-\frac{\|z\|^2}{2}}\|_r\leq \frac{1}{2^{\frac{n}{2r}}\pi^\frac{n}{2}}\|f\|_p\|g\|_q, \ws\ws \frac{1}{r}=\frac{1}{p}-\frac{1}{q'}, \ws\ws 1\leq p\leq \infty, \ws\ws p\leq q'.$$\end{statement}
We have to remark here that in case of $r=p=2$, $q=1$ our constant is greater than in the theorem cited above.
 
Let us see an example now, let $w_i(\zeta_i)=e^{-\zeta_i^2}$, $i=1,\dots,n$. Then it is easy to see that
$$w_{z^0}(\zeta)=\prod_{i=1}^n\frac{1}{\sqrt{\pi}}e^{-2\zeta_i^2-(z^0_i)^2+2\sqrt{2}z^0_i\zeta_i},$$
$$p_{z^0,\alpha}(\zeta)=\prod_{i=1}^n(2\pi)^{\frac{1}{4}}h_{k_i}(\sqrt{2}\zeta_i-z^0_i),$$
where $h_{k_i}$ are the normalized Hermite polynomials.
$$e^{\langle z^0,z\rangle}_{w_{z^0}}=\left(\frac{2}{3}\right)^{\frac{n}{2}}e^{-\frac{\|z^0\|^2+\|z\|^2}{6}+\frac{2}{3}\langle z^0,z\rangle}.$$
$$\varphi_{z^0,\alpha}(z)=K^{\alpha}_{w_{z^0}}\left(e^{\langle z^0,z\rangle}\right)(z)=p_{z^0,\alpha}(L)\left(\left(\frac{2}{3}\right)^{\frac{n}{2}}e^{-\frac{\|z^0\|^2+\|z\|^2}{6}+\frac{2}{3}\langle z^0,z\rangle}\right),$$
and 
$$e^{-\frac{\|z\|^2+\|\zeta\|^2-\|z^0\|^2}{2}+\sqrt{2}\langle z,\zeta\rangle}=\sum_{\alpha}\varphi_{z^0,\alpha}(z)p_{z^0,\alpha}(\zeta)e^{-\|\zeta\|^2+\sqrt{2}\langle z^0,\zeta\rangle}.$$

\medskip

\subsection{Poisson wavelet transform}

This transform was introduced in \cite{k}. The authors pointed out that Poisson wavelet transform is a useful tool to examine exothermic continuous stirred
tank reactions. By "waviness" property, the kernel function of a wavelet transform can not be positive. We chose Poisson wavelet transform, because it's kernel has only one sign change. Finite many sign changes can be handled on the same way. Wavelet systems ensure an expansion of an $f(y)$ and the chromatic method ensure an expansion of $\left(W_{\Psi_n}f\right)(a,b)$ which is similar to the wavelet expansion in some sense.

For each positive integer $n$ the Poisson wavelet $\Psi_n(y)$ on the positive axis is defined by $\Psi_n(y)=-q_n'(y)$, where $q_n(y)=\frac{y^n}{n!}e^{-y}$, that is
$$\Psi_n(y)=\left\{\begin{array}{ll}\frac{y-n}{n!}y^{n-1}e^{-y} \ws \mbox{for}\ws y\geq 0 \\ 0 \ws \mbox{for}\ws y<0\end{array}\right.$$
It is related to the Poisson distribution as it follows. Let $X$ be a discrete random variable having the Poisson distribution with parameter $y$. Then $\mathrm{Prob}(X=n)=q_n(y)$.
The members of this wavelet family have the "waviness" property:
$$\int_{\mathbb{R}}\Psi_n(y)dy=0,$$
and the admissibility constant associated with $\Psi_n(y)$ is 
$$C_{\Psi_n}=\int_{\mathbb{R}}\frac{\left|\mathcal{F}\left(\Psi_n\right)(y)\right|^2}{|y|}dy=\frac{1}{n}.$$
Let us denote the kernel of the wavelet transform by
\begin{equation} u_n(a,b,y)=\frac{1}{\sqrt{|a|}}\Psi_n\left(\frac{y-b}{a}\right),\end{equation}
where $y,b \in \mathbb{R}$ and $a \in \mathbb{R}\setminus \{0\}$. Thus the wavelet transform of a function $f(y)$ (if it exists) is
$$\left(W_{\Psi_n}f\right)(a,b)=\int_{\mathbb{R}}f(y)u_n(a,b,y)dy.$$
Let us define the following differential operator:
\begin{equation}D^{(n)}_{a,b}=\sum_{k=1}^{n+1}\binom{n+1}{k}(-1)^{k+1}a^k\left(a\frac{\partial^k}{\partial a\partial b^{k-1}}+\frac{2k-1}{2}\frac{\partial^{k-1}}{\partial b^{k-1}}\right) + b.\end{equation}
One can check that 
\begin{equation}D^{(n)}_{a,b}u_n(a,b,y)= yu_n(a,b,y).\end{equation}
Indeed, since 
$$\Psi_n(y)+\Psi'_n(y)=\Psi_{n-1}(y), \ws n=2,3, \dots,$$
we have
$$\sum_{k=0}^{n-1}\binom{n-1}{k}\Psi_n^{(k)}(y)=\Psi_1(y),$$
and taking into consideration that 
$$\Psi''_1(y)+2\Psi'_1(y)+\Psi_1(y)=0,$$
we have
\begin{equation}\sum_{k=0}^{n+1}\binom{n+1}{k}\Psi_n^{(k)}(y)=0.\end{equation}
Moreover if $b \in \mathbb{R}$ and $a \in \mathbb{R}\setminus \{0\}$, 
\begin{equation}\frac{\partial^k u_n(a,b,y)}{\partial b^k}=\frac{(-\mathrm{sgn}a)^k}{|a|^{\frac{2k+1}{2}}}\left(\Psi_n^{(k)}\right)\left(\frac{y-b}{a}\right),\end{equation}
and
\begin{equation}\frac{\partial^{k+1} u_n(a,b,y)}{\partial a\partial b^k}=\frac{(-\mathrm{sgn}a)^{k+1}}{|a|^{\frac{2k+3}{2}}}\frac{2k+1}{2}\left(\Psi_n^{(k)}\right)\left(\frac{y-b}{a}\right)$$$$-\frac{(-\mathrm{sgn}a)^k}{|a|^{\frac{2k+5}{2}}}(y-b)\left(\Psi_n^{(k+1)}\right)\left(\frac{y-b}{a}\right).\end{equation}
So substituting (52) into (53) and (54) we get (51), which ensures a chromatic derivation of the wavelet transforms of a function. As above, if $p$ is any polynomial of real variable, $p\left(D^{(n)}_{a,b}\right)u_n(a,b,y)= p(y)u_n(a,b,y)$. Now we turn to the definition of the chromatic derivation of a wavelet transform with respect to a point and weights. Let $(a_0,b_0)$ be a fixed point ($b_0 \in \mathbb{R}$, $a_0 \in \mathbb{R}\setminus \{0\}$). $u_n(a_0,b_0,y)$ has one sign change at $a_0n+b_0$. We divide the domain of the integration to two parts: $I^{(1)}:=I^{(1)}_{(a_0,b_0)}$ and $I^{(2)}:=I^{(2)}_{(a_0,b_0)}$. If $a_0>0$ then $I^{(1)}=(b_0,b_0+a_0n)$ and $I^{(2)}=(b_0+a_0n, \infty)$, if $a_0<0$, then $II^{(1)}=(b_0+a_0n,b_0)$ and $I^{(2)}=(-\infty,b_0+a_0n)$. Further on let $a_0>0$, say. Let $f$ be a function on  $\mathbb{R}$. $f=f^{(1)}+f^{(2)}$, where
$$f^{(1)}(y)=f^{(1)}_{(a_0,b_0)}=\left\{\begin{array}{ll}f(y),\ws\mbox{if}\ws y\in (b_0,b_0+a_0n),\\ 0 \ws\ws \mbox{otherwise}\end{array}\right.$$
and
$$f^{(2)}(y)=f^{(2)}_{(a_0,b_0)}=\left\{\begin{array}{ll}f(y),\ws\mbox{if}\ws y\in (b_0+a_0n,\infty),\\ 0 \ws\ws \mbox{otherwise}\end{array}\right..$$
For $i=1,2$ let $w^{(i)}$ be positive weights with  $\mathrm{supp} w^{(i)}\subset \mathrm{clos} I_i$, such that $w^{(i)}_{n,(a_0,b_0)}$ has finite moments, where
$$\sqrt{w^{(i)}_{n,(a_0,b_0)}}:=\sqrt{w^{(i)}}(-1)^iu_n^{(i)}(a_0,b_0,y).$$
Let $\{p_m^{(i)}\}_m:=\{p^{(i)}_{m,n,(a_0,b_0)}\}_m$ the systems of orthogonal polynomials on $I^{(i)}$ with respect to the weights $w^{(i)}_{n,(a_0,b_0)}$ ($i=1,2$). We define for a suitable $f$
$$K^m_{w^{(i)}_{n,(a_0,b_0)}}\left(W_{\Psi_n}f\right)(a,b):=p_m^{(i)}\left(D^{(n)}_{a,b}\right)\int_{I^{(i)}}f(y)\sqrt{w^{(i)}(y)}u_n(a,b,y)dy$$ $$=\int_{I^{(i)}}f(y)p_m^{(i)}(y)\sqrt{w^{(i)}(y)}u_n(a,b,y)dy.$$
Thus let the $m^{th}$ chromatic derivative of $W_{\Psi_n}f$ with respect to $(a_0,b_0)$ and $w^{(1)},w^{(2)}$
$$K^m_{w^{(1)},w^{(2)}(a_0,b_0)}\left(W_{\Psi_n}f\right)(a,b)$$ \begin{equation}:=K^m_{w^{(1)}_{n,(a_0,b_0)}}\left(W_{\Psi_n}f\right)(a,b)+K^m_{w^{(2)}_{n,(a_0,b_0)}}\left(W_{\Psi_n}f\right)(a,b).\end{equation}
Let $f\in L^2(\mathbb{R})$. Let $g^{(i)}(y)=g^{(i)}_{n, w^{(i)},(a_0,b_0)}=\frac{f^{(i)}}{\sqrt{w^{(i)}_{n,(a_0,b_0)}}}$. 
$$g^{(i)}\sim \sum_{m=0}^{\infty}c_m^{(i)}p_m^{(i)}, \ws \mbox{where} \ws c_m^{(i)}=\int_{I^{(i)}}g^{(i)}p_m^{(i)}w^{(i)}_{n,(a_0,b_0)}.$$
So $f$ has Fourier expansions in both intervals
$$f= \sum_{m=0}^{\infty}c_m^{(1)}p_m^{(1)}\sqrt{w^{(1)}_{n,(a_0,b_0)}}+\sum_{m=0}^{\infty}c_m^{(2)}p_m^{(2)}\sqrt{w^{(2)}_{n,(a_0,b_0)}}.$$
Let us observe that
$$c_m^{(1)}=-K^m_{w^{(1)}_{n,(a_0,b_0)}}\left(W_{\Psi_n}f^{(1)}\right)(a_0,b_0), \ws \mbox{and} \ws c_m^{(2)}=K^m_{w^{(2)}_{n,(a_0,b_0)}}\left(W_{\Psi_n}f^{(2)}\right)(a_0,b_0).$$
Let
$$\varphi(a,b):=\varphi_{(a_0,b_0),n}(a,b)=\left(W_{\Psi_n}\right)\left(u_n(a_0,b_0,y)\right)(a,b)$$
and
$$\varphi_m^{(i)}(a,b):=\varphi^{(i)}_{m,(a_0,b_0),n}(a,b)=\left(W_{\Psi_n}\right)\left(p_m^{(i)}(y)\sqrt{w^{(i)}_{n,(a_0,b_0)}(y)}\right)(a,b).$$
Thus
$$K^m_{w^{(i)}_{n,(a_0,b_0)}}\left(\varphi\right)(a,b)=(-1)^i\varphi_m^{(i)}(a,b).$$
By Parseval's formula
$$n\int_{\mathbb{R}}\int_{\mathbb{R}}\varphi_l^{(i)}(a,b)\varphi_m^{(i)}(a,b)db\frac{da}{a^2}=\delta_{l,m},$$
and
$$c_m^{(i)}=\int_{\mathbb{R}}f^{(i)}(y)p_m^{(i)}(y)\sqrt{w^{(i)}_{n,(a_0,b_0)}}(y)dy=n\int_{\mathbb{R}}\int_{\mathbb{R}}\left(W_{\Psi_n}f^{(i)}\right)(a,b)\varphi^{(i)}_{m}(a,b)db\frac{da}{a^2}.$$
That is $\left(W_{\Psi_n}f^{(i)}\right)(a,b)$ has an orthogonal expansion with respect to $\{\varphi_m^{(i)}(a,b)\}_m$,
$$\left(W_{\Psi_n}f^{(i)}\right)(a,b)\sim \sum_{m=0}^{\infty}c_m^{(i)}\varphi_m^{(i)}(a,b)$$ 
$$=\sum_{m=0}^{\infty}(-1)^iK^m_{w^{(i)}_{n,(a_0,b_0)}}\left(W_{\Psi_n}f^{(i)}\right)(a_0,b_0)(-1)^iK^m_{w^{(i)}_{n,(a_0,b_0)}}\left(\varphi\right)(a,b),$$
that is 
$$\left(W_{\Psi_n}f\right)(a,b)\sim \sum_{m=0}^{\infty}K^m_{w^{(1)}_{n,(a_0,b_0)}}\left(W_{\Psi_n}f\right)(a_0,b_0)K^m_{w^{(1)}_{n,(a_0,b_0)}}\left(\varphi\right)(a,b)$$
\begin{equation}+\sum_{m=0}^{\infty}K^m_{w^{(2)}_{n,(a_0,b_0)}}\left(W_{\Psi_n}f\right)(a_0,b_0)K^m_{w^{(2)}_{n,(a_0,b_0)}}\left(\varphi\right)(a,b)\end{equation}
Since $u_n(a,b,y)$ is in $L^2(\mathbb{R})$ as a function of $y$, it is easy to see
\begin{equation}u_n(a,b,\cdot)= \sum_{m=0}^{\infty}\varphi_m^{(1)}(a,b)p_m^{(1)}\sqrt{w^{(1)}_{n,(a_0,b_0)}}+\sum_{m=0}^{\infty}\varphi_m^{(2)}(a,b)p_m^{(2)}\sqrt{w^{(2)}_{n,(a_0,b_0)}}.\end{equation}
Let us denote by
$$S_N(a,b)=\sum_{m=0}^{N}\left(c_m^{(1)}\varphi_m^{(1)}(a,b)+c_m^{(2)}\varphi_m^{(2)}(a,b)\right)=S_N^{(1)}(a,b)+S_N^{(2)}(a,b),$$
and by
$$s_N(y)=\sum_{m=0}^{N}\left(c_m^{(1)}p_m^{(1)}\sqrt{w^{(1)}_{n,(a_0,b_0)}}+c_m^{(2)}p_m^{(2)}\sqrt{w^{(2)}_{n,(a_0,b_0)}}\right)=s_N^{(1)}(y)+s_N^{(1)}(y),$$
Let
$$\|h(a,b)\|_2==\left(n\int_{\mathbb{R}}\int_{\mathbb{R}}\left|h(a,b)\right|^2db\frac{da}{a^2}\right)^{\frac{1}{2}}.$$
Following the chain of ideas of \cite{z1}, we have

\medskip 

\begin{theorem} Let $f\in L^2(\mathbb{R})$.
$$\lim_{N\to\infty}\left\|\left(W_{\Psi_n}f\right)(a,b)-S_N(a,b)\right\|_2=0.$$
If  $a\neq 0$ and $b \in \mathbb{R}$, then
$$\left|\left(W_{\Psi_n}f\right)(a,b)-S_N(a,b)\right|\longrightarrow 0, \ws\ws (N \to \infty) $$
on $\mathbb{R}\setminus\{0\}\times  \mathbb{R}$.
\end{theorem}

\medskip 

\proof By Parseval's formula
$$\left(n\int_{\mathbb{R}}\int_{\mathbb{R}}\left|\left(W_{\Psi_n}f\right)(a,b)-S_N(a,b)\right|^2db\frac{da}{a^2}\right)^{\frac{1}{2}}=\|f-s_N\|_2$$ $$\leq \|f^{(1)}-s_N^{(1)}\|_2+\|f^{(2)}-s_N^{(2)}\|_2,$$
which tends to zero by Bessel's inequality. Moreover by (56)
$$\sum_{m=0}^{\infty}\left(\varphi_m^{(i)}(a,b)\right)^2=\|u_n^{(i)}(a,b,y)\|_2\leq \|u_n(a,b,y)\|_2\leq C(n),$$
and so
$$\left|\left(W_{\Psi_n}f^{(i)}\right)(a,b)-S_N^{(i)}(a,b)\right|^2\leq C(n)\sum_{m=N+1}^{\infty}\left|c_m^{(i)}\right|^2.$$
Thus the triangle inequality ensures the second statement of the theorem.

\medskip

\noindent \small{Department of Analysis, \newline
Budapest University of Technology and Economics}\newline
\small{ g.horvath.agota@renyi.mta.hu}

\end{document}